\documentclass[12pt,a4paper]{amsart}
\usepackage{graphicx,amssymb}
\usepackage[usenames]{color}
\input xy
\xyoption{all}
\usepackage[all]{xy}
\usepackage{hyperref}

\newtheorem{thm}{Theorem}[section]
\newtheorem{cor}[thm]{Corollary}
\newtheorem{lem}[thm]{Lemma}
\newtheorem{prop}[thm]{Proposition}
\theoremstyle{definition}

\newtheorem{rem}[thm]{Remark}

\numberwithin{equation}{section}

\newcommand{\QQ}{\mathbb Q}

\newcommand{\ZZ}{\mathbb Z}
\newcommand{\CC}{\mathbb C}
\newcommand{\PP}{\mathbb P}

\newcommand{\lra}{\longrightarrow}

\newcommand{\ra}{\rightarrow}

\newcommand{\cH}{\mathcal{H}}

 \DeclareMathOperator{\Nm}{{Nm}}

 \DeclareMathOperator{\Ima}{{Im}}

 \DeclareMathOperator{\tr}{tr}

\DeclareMathOperator{\End}{{End}}

\DeclareMathOperator{\Irr}{{Irr}}
\begin{document}

\title[ ]{ Decomposable Jacobians}
\author{Angel Carocca, Herbert Lange and Rub\'i E. Rodr\'iguez}
\address{Departamento de Matem\'atica y Estad\'istica, Universidad de La Frontera, Avenida Francisco Salazar 01145, Casilla 54-D, Temuco, Chile.}
\email{angel.carocca@ufrontera.cl}
\address{Department Mathematik, Universit\"at Erlangen, Cauerstrasse 11, 91058 Erlangen, Germany.}
\email{lange@math.fau.de}
\address{Departamento de Matem\'atica y Estad\'istica, Universidad de La Frontera, Avenida Francisco Salazar 01145, Casilla 54-D, Temuco, Chile.}
\email{rubi.rodriguez@ufrontera.cl}

\thanks{The authors were partially supported by Grants Fondecyt 1190991 and CONICYT PAI Atracci\'on de Capital Humano Avanzado del Extranjero PAI80160004}
\begin{abstract}
In this paper we give examples of smooth projective curves whose Jacobians are isogenus to a product of an arbitrarily high 
number of Jacobians.
\end{abstract}

\maketitle

\section{Introduction}

In \cite{es},  Ekedahl and Serre gave examples of curves $X$ whose Jacobian is completely decomposable, i.e.
 isogenous to a product of elliptic curves.  The 
highest genus of their examples is  $g(X) =1297$. They asked among other things whether the genus of a curve with a
completely decomposable Jacobian is bounded above. Although in the meantime  many other examples of such Jacobians 
have been given (see e.g. \cite{pr}), no example of genus bigger than 1297 seems to be known.  

In this paper we consider an easier question, namely: can a Jacobian be isogenous to the product of arbitrary many Jacobians 
of the same genus (not necessarily equal to one)? In \cite{clr1} we gave examples of Jacobians which are isogenous to an 
arbitrary number of Prym varieties of the same dimension. The main result of this paper is the following theorem (see Corollary 
\ref{c4.6}).

\begin{thm}
Given any positive integer $N$, there exists a smooth projective curve $X$ whose Jacobian is isogenous to the product of
 $m \ge N$ Jacobian varieties of the same dimension.
\end{thm}

To be more precise, in Corollary \ref{c4.5}, for any positive integer $n \geq 2$ 
and any   integer $t \ge 2n$  we give examples of curves of genus $m(t-3)$ whose Jacobians are isogenous to the product of  $m = \frac{1}{3}(2^{2n} -1)$ Jacobians of curves of genus $ \; t - 3.$
The idea is to consider curves $X$ of genus  $2^{2n}(t-3) +1$ with an action of the group 
$G = N \rtimes B \simeq \ZZ_2^{2n} \rtimes \ZZ_3$ of signature $(0;\underbrace{3,\dots,3}_{t})$.  Then the quotient
curves  $X_B = X/B$
will have the asserted properties.\\

The $m$ Jacobians are Jacobians of subcovers of the curve $X, $ of genus $t - 3,$ and they are 
in fact isomorphic to Prym varieties of  \'etale double covers of other subcovers. Using these isomorphisms gives a second proof of
the theorem, which actually gives a bit more, namely an estimate of the degree of the isogeny.\\

In Section 2 we investigate the representations of the group $G$. In Section 3 we recall some results on curves with $G$-action. In Section 4 we study the diagram of subcovers of the curve $X$. Section 5 contains the first proof of the above mentioned
theorem. In Section 6 we give the proof using the Prym varieties and the trigonal construction. Finally in the last section we 
compute the group algebra decomposition of the Jacobian $JX$ of the curve $X$.

\section{The group $G$}

Consider the group $ V_4 \rtimes B  \simeq A_4$ the alternating group of degree 4 with 
$V_4 = \langle a_1, a_2 \rangle$ the Klein group and $B = \langle b \rangle \simeq \ZZ_3$, and for any positive integer $n$  let
$$
G := N \rtimes B :=  \left(\langle a_1, a_2 \rangle \times  \langle a_3, a_4 \rangle \times \cdots \times  \langle a_{2n-1}, a_{2n} \rangle \right) \rtimes 
\langle b \rangle
$$
where the subgroups $\langle a_j,a_{j+1} \rangle \simeq V_4$ for each $j$ and $\langle a_j,a_{j+1} \rangle \rtimes  \langle b \rangle \simeq A_4$ for $j = 1,3,\dots 2n-1$. Obviously $G$ is of order
$$
|G| = 2^{2n}3.
$$
Also $ \; B \; $ is a $3$-Sylow subgroup of $ \; G \; $ and $ \; B = N_G(B) = C_G(B).$\\

In this section we will show that all subgroups $M$ of $G$ of index $4$ are maximal in $G$ (there is no subgroup strictly between it and $G$); if we consider only such $M$ containing $B$ (there is exactly one in each conjugacy class of subgroups of index four), they correspond bijectively to the subgroups $U$ of index four in $N$ that are normal in $G$, and also to the set of irreducible representations of $G$ of degree three. An example of such a subgroup is $ \; M_1 = \left(  \langle a_3, a_4 \rangle \times \cdots \times  \langle a_{2n-1}, a_{2n} \rangle \right) \rtimes 
\langle b \rangle$ with $U_1 = \langle a_3, a_4 \rangle \times \cdots \times  \langle a_{2n-1}, a_{2n} \rangle$. We need some preliminaries.

\begin{lem} \label{l1.1}
 $ \; G \; $ has no subgroups of index $ \; 2. \; $
\end{lem}

\begin{proof} Suppose $ \; R \leq G \; $ such that $ \; \vert G : R \vert = 2. \; $ Then $ \; R \unlhd G \; $ and by the Sylow Theorem $ \; B \leq R. \; $ In this way,  we have $ \; G = NB = NR. \;$ This implies
$$ Q := N \cap R \unlhd G 
$$ 
and $ \; \; \vert N : Q \vert = 2. \; $

Hence $ \; N = \langle n_1 \rangle  Q \; $ for some $n_1 \in N$. Since $ \; Q \unlhd G, \; $ we have 
$$    (n_1)^b := bn_1b^{-1} = n_2 \: \not \in   \;  Q \quad  \mbox{and} \quad (n_1)^{b^2} = n_1n_2   \not \in   Q.  
$$ 
The last equation follows from the fact that $\langle  n_1,n_2 \rangle \simeq V_4$ and $ (n_1)^{b^2} \neq  n_1,n_2$.
Then $  \langle n_1, n_2 \rangle \cap Q = \{1 \}$  and hence
$$  \; 2 = \vert N : Q \vert = \displaystyle\frac{\vert \langle n_1, n_2 \rangle \vert}{\vert \langle n_1, n_2 \rangle \cap Q \vert } = 4, 
$$ 
a contradiction.
\end{proof}

Now let $ \; M \leq G \; $  be any subgroup such that $ \; \vert G : M \vert = 4. \; $  As we saw above, there are such subgroups. According to Lemma \ref{l1.1}, $M$ is maximal in $G$. Moreover, we have,

\begin{lem}
$ \; G = NM \; \; $ ; $ \; U = N \cap M \unlhd G \; $  and $ \; \vert N : U \vert = 4.$
\end{lem}

\begin{proof} Since $ \; M \; $ is a maximal subgroup of $ \; G \; $ and $ \; N \unlhd G \; $ we have  $ \; G = NM. \; $ \\
Since $ \; N \unlhd G \; $ we have  that $ \;  U = N \cap M \unlhd M. \; $ Also, $ \; U \unlhd N \; $ since $ \; N \; $ is abelian. 
Therefore, $ \; U \unlhd NM = G. \; \; $ Also 
 $ \; 4 = \vert G : M \vert = \vert N : U \vert.$
\end{proof}

According to a result of elementary number theory, the number 
$$
m := \frac{2^{2n} -1}{3}
$$
is an integer for all positive integers $n$.
\\

\begin{lem}
There are complex  irreducible representations $\theta_1, \dots ,\theta_m$ of degree $3$ and complex irreducible representations 
$\chi_0, \chi_1, \chi_2$ of degree $1$ of $G$ ($\chi_0$ the trivial character). These are all the irreducible complex representations of $G$.  
\end{lem}

\begin{proof}
For any  irreducible character $\psi$ of $N$ and any element $g \in G$  we denote by $\psi^g$
the conjugate character of $N$ defined by $\psi^g(n) := \psi(gng^{-1})$, 
The stabilizer of the trivial character  $\psi_0 $ of $N$
is the group $G$, whereas the stabilizer of any non-trivial character of $N$ is trivial. Hence there are exactly $ 1 + m$ orbits
for the action of $G$ on  the set of all irreducible characters of $N$. Let $\psi_0,  \psi_1, \dots, \psi_m $ be a system of representatives of these 
orbits and $ \; \nu_0, \nu_1, \nu_2 \; $ the irreducible characters of $ B. \; $  Then, according to  \cite[Proposition 25]{s} the irreducible representations  of $ \; G \; $ are 
$ \; \chi_0 = \psi_0 \otimes \nu_0, \; \chi_1 = \psi_0 \otimes \nu_1, \; \chi_2 = \psi_0 \otimes \nu_2 \; $ and the induced representation $ \; \theta_i = (\psi_i \otimes \nu_0)^G, \;  $ for each $ \; i = 1, \dots m.$
\end{proof}

This implies immediately,

\begin{cor} \label{c2.4}
The irreducible rational representations of $G$ are exactly the trivial representation $\chi_0$, the representations $\theta_j$ of degree $3$
for $j = 1, \dots , m$, and the representation $\psi := \chi_1 \oplus \chi_2$ of degree $2$.
\end{cor}

\begin{lem}  Let $\psi$ be a non-trivial irreducible character of $N$ and $\theta:= (\psi \otimes \nu_0)^G $ the corresponding irreducible $3$-dimensional representation of $G$. Then
\begin{enumerate}
\item[(i)] $  \vert N : \ker(\psi) \vert = 2; \; $
\item[(ii)] $ \; U := \ker(\psi) \cap \ker(\psi^b) = \ker(\psi)  \cap \ker(\psi\psi^{b}) \; \; $ ;  $ \; \vert N : U \vert = 4 \; \; $ and $ \; \; \ker(\theta) = U \unlhd G.$ 
\end{enumerate}
\end{lem}

\begin{proof} (i):  Since $ \; N /\ker(\psi) \; $ is isomorphic to a finite cyclic subgroup of $ \; {\mathbb C}^{\ast} \; $ and $ \; N \cong ({\mathbb Z}_2)^{2n} \; $ we have $ \; \vert N : \ker(\psi) \vert = 2. \; $
\\
(ii): Since $ \; \vert N : \ker(\psi) \vert = 2 \; $  and $ \; \psi \ne \psi^b \; $ we have $ \; N = \ker(\psi) \ker(\psi^b). $ \\
Hence $  \; 2 = \vert N : \ker(\psi) \vert = \vert \ker(\psi^b) : \ker(\psi) \cap \ker(\psi^b) \vert. \; $ Therefore
 $ \; \vert N : U \vert = 4.$

We have $ \; \ker(\psi) \cap \ker(\psi^b) \leq \ker(\psi \psi^b). \; $ 
So 
$$  U = \ker(\psi) \cap \ker(\psi^b)  = \ker(\psi) \cap \ker(\psi^b) \cap \ker(\psi \psi^b) = \cap_{i=0}^2 \ker (\psi^i).
$$ 
It is well known that $\ker(\theta) = \bigcap_{g \in G}\ker(\psi^g)$. Hence
$$\ker(\theta) = \bigcap_{g \in G}\ker(\psi^g) = \bigcap_{n \in N, \; i = 0,1,2 }\ker(\psi^{n{b^i}}) =  \bigcap_{i=0}^{2}\ker(\psi^{{b^i}})  = U. $$
\end{proof}

Conversely, we have

\begin{lem} \label{l2.4}
 Let $ \; U \unlhd G \; $ such that $ \; \vert N : U \vert = 4. \; $ Then
\begin{enumerate}
\item[(i)]  There is a non-trivial character$ \; \psi \; $ of $ \; N \; $ such that $ \; \theta = (\psi \otimes \nu_0)^G \; $ is an irreducible representation of $G$ with $ \; U = \ker(\theta);$ 

\item[(ii)] $ \; M = UB \leq G \; $ and $ \; \vert G : M \vert = 4.$
\end{enumerate}
\end{lem}

\begin{proof} 
(i): We have that $ \; N/U \; $ is isomorphic to the Klein group of order $ \; 4. \; $ Consider $ \; L \leq N \; $ such that $ \; U \leq L \; \; $ and $ \; \vert N : L \vert = 2. \; $ Then $ \; U \leq L^{b^i} \; $ since $ \; U \unlhd G.$

There is a non-trivial character $ \; \psi \; $ of $ \; N \; $ such that $ \; L = \ker(\psi). \; $ In this way, for the induced representation $ \; \theta = (\psi \otimes \nu_0)^G, \; $ we have by the previous lemma that $ \; U = \ker(\theta).$

(ii): Since $ \;  U \unlhd G \; $ with $[N:U] =4$, we have $ \; M = UB \leq G$ and $|G:M| =4$.
\end{proof}

Combining these lemmas  we get
\begin{prop} \label{p2.5}
\begin{enumerate} 
\item[(i)] There are canonical bijections between the following sets
\begin{itemize}
\item $\{ U \unlhd G \; \; / \; \; \vert N : U \vert = 4 \} ;$
\item $\{ M \leq G \; \; / \; \; \vert G : M \vert = 4 \ \textup{ and }  \ B \leq M \} ; $ 
\item $ \{ \theta \in \Irr(G) \; \; / \; \;  \deg(\theta) = 3 \}.$
\end{itemize}
\item[(ii)]  $$ 
\vert \{ \theta \in \Irr(G) \; \; / \; \;  \deg(\theta) = 3 \} \vert = m = \displaystyle{\frac{2^{2n} - 1}{3}}
$$
\end{enumerate}
\end{prop}

\begin{proof}
(i): The bijections are given by
$$
U = M \cap N  \longleftrightarrow M = UB \longleftrightarrow U = \ker(\theta).
$$	
If $ \; \rho_M \; $ denotes the  representation of $G$ induced by the trivial representation of $ \; M, \; $ then for each such $M$,
$$
\rho_M = \chi_0 \oplus \theta.
$$
So $\theta$ is of degree 3.

(ii): There are $2^{2n} - 1$ non-trivial irreducible characters $\psi$ of $N$. The group $B$ acts non-trivially on them. So Lemma 
\ref{l2.4} and part (i) of the proposition imply the assertion.
\end{proof}

\begin{rem} \label{r2.6}
Let $ \; U \unlhd G \; $ such that $ \; \vert N : U \vert = 4,   \; $   $ \; M = UB  \; $ (the corresponding maximal subgroup), and $U = \ker(\theta)$ as in the Proposition \ref{p2.5}.

Consider $ \; L \leq N \; $ such that $ \; U \leq L  $ and $ \; \vert N : L \vert = 2.$ Observe that for each $U$ there are three such $L$, forming a conjugacy class of subgroups of $G$.
Also, consider $ \; \psi  $, the irreducible character of $ \; N \; $ with $ \; \ker(\psi) = L. \; $ Then
\begin{itemize}
\item $ \; U \leq \ker(\psi^{b^i}) =: L_i \; $ for $ \; i = 0,1,2.\;$ So $L_0 = L$;
\item $ \; \theta = (\psi^{b^i} \otimes \nu_0)^G \; $ for $ \; i = 0,1,2 $;
\item $ \; \rho_{N} = \chi_0 \oplus \chi_1 \oplus \chi_2 $;
 \item $ \; \rho_{L_i} =  \ \rho_N \oplus \theta \; $ for $ \; i = 0,1,2.$;
 \item $ \; \rho_M = \chi_0 \oplus \theta .$
\end{itemize}
\end{rem}

\section{Action of Hecke algebras of an abelian variety}

Let $G$ be any finite group acting on an abelian variety $A$ over the field of complex numbers and let $H \le G$ be a subgroup. 
In this section we recall a result of \cite{clr1} together with the notation needed for it.

The element 
$$
p_H := \frac{1}{|H|} \sum_{h \in H} h
$$
is an idempotent of the group algebra $\QQ[G]$. The {\it Hecke algebra} for $H$ in $G$ is defined to be the subalgebra
$$
\cH_H := p_H\QQ[G]p_H = \QQ[H\backslash G/H] 
$$
of $\QQ[G]$. The action of $G$ on $A$ induces an algebra homomorphism
\begin{equation} \label{e3.1}
\QQ[G] \ra \End_\QQ(A)  
\end{equation}
in a natural way.  We denote the elements of $\QQ[G]$ and their images by the same letter. For any element $\alpha \in \QQ[G]$ we
define its image in $A$ by 
$$
\Ima(\alpha) :=\Ima(n\alpha) \subset A
$$
where $n$ is any positive integer such that $n\alpha$ is in $\End(A)$. It is an abelian subvariety which does not depend on the 
chosen integer $n$. Consider the abelian subvariety of $A$ given by
$$
A_H := \Ima(p_H).
$$
Restricting \eqref{e3.1} to $\cH_H$ gives an algebra homomorphism $\cH_H \ra \End_\QQ(A_H).$

Let $\{W_1, \dots, W_r\}$ denote the irreducible rational representations of $G$.
To any $W_i$ there corresponds an
 irreducible complex representation $V_i$, uniquely determined up to an element of the Galois group of $K_i$ over $\QQ$,
where $K_i$ is the field obtained by adjoining to $\QQ$ the values of the character $\chi_{V_i}$  of $V_i$.

To each $W_i$ we can associate a central idempotent $e_{W_i}$ of $\QQ[G]$ by
$$
e_{W_i} = \frac{\dim_\CC(V_i)}{|G|} \sum_{g \in G}  \tr_{K_i/\QQ}(\chi_{V_i}(g^{-1}))g \, .
$$
Let $\rho_H$ denote the representation of $G$ induced by the trivial representation of $H$. It decomposes as
\begin{equation} \label{e2.1}
\rho_H \simeq \sum_{i=1}^r a_i W_i,
\end{equation}
with $a_i = \frac{1}{s_i} \dim_\CC(V_i^H)$ and $s_i$ the Schur index of $V_i$.  Renumbering if necessary, let $\{W_1, \dots, W_u\}$ denote the set of all 
irreducible  rational representations of $G$ such that $a_i \neq 0$. Then there is a bijection  from this set to the set 
$\{ \widetilde W_1, \dots, \widetilde W_u\}$ of all irreducible rational representations of the algebra $\cH_H$.
An analogous statement holds for the complex  irreducible representations of $G$ and of $\cH_H$. Let $\widetilde V_i$ denote the  
representation of  $\cH_H$ associated to the complex irreducible representation $V_i$ of $G$ and Galois associated to $ \widetilde{W}_i. $
The dimension of $\widetilde W_i$ is given by
\begin{equation} \label{eq2.2}
\dim_\QQ(\widetilde W_i) = \dim_\QQ (W_i^H) = [L_i:\QQ] \dim_\CC({V}^H_i).
\end{equation}
where $L_i$ denotes the field of definition of the representation $V_i$. Recall that $s_i = [L_i:K_i]$ is the Schur index of $V_i$.

For $i = 1, \dots,u$ consider the central idempotents of $\cH_H$ given by
$$
f_{H,\widetilde W_i} := p_H e_{W_i} = e_{W_i} p_H.
$$ 
Then $p_H$ decomposes as 
\begin{equation} \label{e3.4}
p_H = \sum_{i=1}^u f_{H,\widetilde W_i} \, .
\end{equation}
Defining for $i=1, \dots , u$ the abelian subvarieties 
$$
A_{H, \widetilde W_i} := \Ima(f_{H,\widetilde W_i}),
$$
one obtains the following isogeny decomposition of $A_H$,  given by the addition map
\begin{equation} \label{e2.2}
+:   A_{H, \widetilde W_1} \times A_{H, \widetilde W_2} \times \cdots \times A_{H, \widetilde W_u}  \ra A_H.
\end{equation}
It is uniquely determined by $H$ and the action of $G$ and  called the {\it isotypical decomposition of} $A_H$.
So in order to describe the abelian subvariety $A_H$, it suffices to describe the $A_{H,\widetilde W_i}$.
This is done by \cite[Theorem 4.3]{clr1}. For it we need one more notation. Consider the decomposition of $G$ into double cosets 
of $H$ in $G$,
$$
G = H_1 \cup H_2 \cup \cdots \cup H_s
$$
with $H_1 = H$. Then a basis for the Hecke algebra $\cH_H$ is given by the elements
$$
q_j := \frac{1}{|H|} \sum_{h_j \in H_j} h_j 
$$
for $j = 1, \dots,s$. If $\chi_{V_i}$ denotes the character of the representation $V_i$, then \cite[Theorem 4.3]{clr1} says,
\begin{thm} \label{thm3.1}
Suppose $\dim V_i^H =1$ and $K_i = \QQ$. Then
$$
A_{H, \widetilde W_i} = \{ z \in A \;|\; q_j(z) = \chi_{V_i}(q_j)z \; \mbox{for} \; j = 1, \dots,s \}_0
$$
where the index $0$ means the connected component containing $0$.
\end{thm}

\section{The diagram of subcovers of $X$}

We start this section recalling some basic properties on groups action on  smooth projective curves.
Let $\; X \;$ be a smooth projective curve of genus $ \;g ,\;  G$ a finite group acting
on $\; X   \; $  and  $\; X \rightarrow X_{G}\; $ be the quotient projection.
This cover may be partially characterized
by a vector of numbers $\;(\gamma ; m_{1}, \cdots , m_{r})\;$ where $\;\gamma \; $ is the genus of $\; X_{G}, $ the integer  $ \; 0 \leq r \leq 2 g + 2 \;$ is the number of branch points of the cover and the integers  $ \; m_{j} \; $ are the orders 
of the cyclic subgroups $ \; G_j \; $ of $ \; G \; $ which fix points on $ \; X$. We call $(\gamma ; m_{1}, \cdots , m_{r})$  the \textit{branching data} of $G$ on $ \; X$. These numbers satisfy the Riemann-Hurwitz formula
\begin{equation}\label{rh}
\frac{2(g -1)}{\vert G\vert }=2(\gamma-1)+\sum_{j=1}^{r}\left(1-\frac{1}{m_{j}}\right) .
\end{equation}
\noindent 
A $(2\gamma+r)-$tuple $\left(a_{1},\cdots, a_{\gamma},b_{1},\cdots,b_{\gamma}, c_{1},\cdots,c_{r}\right)$ of elements of $ \; G \; $ is called a \textit{generating vector of type } $(\gamma;m_{1},\cdots,m_{r})$  if
$$\label{gvector}
 G = \left\langle a_{1},\cdots, a_{\gamma},b_{1},\cdots,b_{\gamma}, c_{1},\cdots,c_{r} \; \; / \;  \prod_{i=1}^{\gamma}[a_{i},b_{i}]\prod_{j=1}^{r}c_{j} = 1 \; ,  \; \vert c_{j}\vert  = m_{j} \;  \mbox {for} \; j =1,...,r \;   \right\rangle\\
$$
where   $[a_i,b_i]=a_ib_ia_i^{-1}b_i^{-1}.$ 
\vspace{2mm}\\
Riemann's Existence Theorem gives us the following theorem   (see \cite[Proposition 2.1]{br})
\begin{thm}
The group $G$ acts on a smooth projective curve $ X \; $ of genus $\; g$ with branching data $(\gamma;m_{1},\cdots,m_{r})$ if and only if $G$ has a generating vector of type $(\gamma;m_{1},\cdots,m_{r})$ satisfying the Riemann-Hurwitz formula \eqref{rh}.
\end{thm}

From now on let  $G = N \rtimes B \simeq \ZZ_2^{2n} \rtimes \ZZ_3$ be the group of Section 2. Recall that
$$G := N \rtimes B :=  \left(\langle a_1, a_2 \rangle \times  \langle a_3, a_4 \rangle \times \cdots \times  \langle a_{2n-1}, a_{2n} \rangle \right) \rtimes 
\langle b \rangle
$$
where for each $j = 1,3, \dots 2n - 1 $ the subgroups $\langle a_j,a_{j+1} \rangle $ are isomorphic to Klein group and the subgroups $\langle a_j,a_{j+1} \rangle \rtimes  \langle b \rangle $ are isomorphic to the alternating group of degree four.
\begin{rem} For $ \; n = 1,  \; $ we have that $ \; G $ is isomorphic to the  alternating group of degree four. It is easy to check that $G$  acts on a smooth projective curve  
$X$ of genus $ \; g = 4(t - 3) + 1 \; $ with branching data $(0;\underbrace{3,\cdots, 3}_{t}) \; $ and $ \; t \geq 4.$\\
In this case, if $ \; 
P(X_{\langle a_1 \rangle} \ra X_N)$ denotes the Prym variety of $X_{\langle a_1 \rangle } \ra X_N $ then, by the trigonal construction (see \cite{rr}), we have
$ P(X_{\langle a_1 \rangle} \ra X_N)$ is isomorphic to $JX_B$  as principally polarized abelian variety.
\end{rem}
\begin{lem}\label{l4.3} For $ \; n > 1 \; $ the group $G$ acts on a smooth projective curve  $X$ of genus $ \; g = 2^{2n}(t - 3) + 1 \; $ with branching data $(0;\underbrace{3,\cdots, 3}_{t}) \; $ for $ \; t \geq 2n.$
\end{lem}
\begin{proof}  We illustrate by giving a generating vector for $ \; t = 2n.$
\vspace{2mm}\\
For $ \; n = 2, \; $ consider the elements $ \; y_1 = a_1b, \; y_2 =b, \; y_3 = a_{3}b^2 \; $ and $ \; y_4 = b^2a_{4}a_1. \; $ \\
Since $ \; b^2a_{3}b = a_{4} \; $ we have $ \; y_1y_2y_3y_4 = a_1bba_{3}b^2 b^2a_{4}a_1 = a_1b^2a_{3}ba_{4}a_1 = a_1a_{4}a_{4}a_1 = 1. \;$ Also, it is clear that
$ \; G = (\langle a_1, a_2 \rangle \times  \langle a_3, a_4 \rangle) \rtimes \ZZ_3 =  \langle y_1, y_2, y_3, y_4 \rangle. \; $ In this way $ \; (y_1,y_2,y_3,y_4) \; $ is a generating vector of type $ \; (0; 3,3,3,3) \; $ for $ \; G.$
\vspace{2mm}\\
For $ \; n = 3, \; $  is it clear that  $ \; ( y_1 = a_1b, \; y_2 =b, \; y_3 = a_{3}b^2,\; y_4 = b^2a_{4}a_1, \; y_5 = a_5b, y_6 = b^2a_5) \; $ is a generating vector of type $ \; (0; 3,3,3,3,3,3) \; $ for 
$ \; G = (\langle a_1, a_2 \rangle \times  \langle a_3, a_4 \rangle \times  \langle a_5, a_6 \rangle) \rtimes \ZZ_3.$
\vspace{2mm}\\
Now, for $ \; n > 3 \; $ we can apply a similar  procedure as the one described  above in both cases $n $ even or odd.

\end{proof}

Let $X$ denote a smooth projective curve with an action of $G$ with signature $(0;\underbrace{3,\dots,3}_t).$  
The action induces an action of $G$ on the Jacobian 
$JX$ of $X$. In this section we will study the subcovers of $X$ which will be important for us.\\

For any subgroup $H$ of $G$
let $X_H$ denote the quotient curve
$$
X_H := X/H.
$$
According to Proposition \ref{p2.5} the subgroup $N$ has exactly $m = \frac{2^{2n}-1}{3}$   subgroups of index  
$ 4 \: $ which are normal in $G$.
Let $U_1, \dots, U_m$ be  these subgroups and $ \; M_1, \dots , M_m \; $ denote the corresponding maximal subgroups of $ \; G.$ For any $U_j$ we choose a subgroup $L_j \le N$ containing $U_j$ with $[N:L_j]=2$. According to Remark \ref{r2.6}, the conjugacy class of $L_j$ is uniquely determined.  With these notations
we have the following diagram of covers of curves.

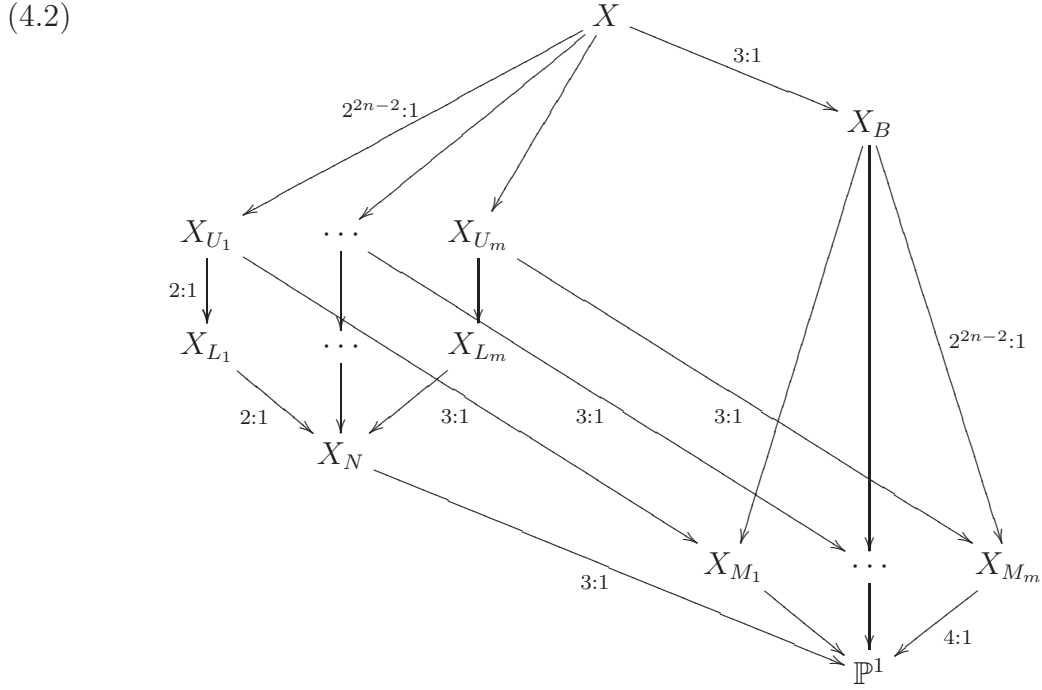
\begin{equation} \label{d3.1}
\xymatrix{
&&& \ar[ddlll]_{2^{2n-2}:1}\ar[ddll]\ar[ddl]X\ar[drr]^{3:1} & & &\\
&&&&&X_B\ar[ddddl]\ar[dddd]\ar[ddddr]^{2^{2n-2}:1} & \\
\ar[d]X_{U_1}\ar[d]_{2:1}\ar[dddrrrr]_{3:1}& \cdots\ar[d] \ar[dddrrrr]_{3:1} & X_{U_m}\ar[d] \ar[dddrrrr]_{3:1} &&&&&\\
X_{L_1}\ar[dr]_{2:1}& \cdots \ar[d]& X_{L_m}\ar[dl] &&&&&\\
 & X_N\ar[ddrrrr]_{3:1}&  &&&&&\\
 & &  &&X_{M_1}\ar[dr]&\cdots \ar[d]&X_{M_m}\ar[dl]^{4:1}\\
  & &&&  &{\mathbb P}^1&\\
}
\end{equation}

\begin{lem} \label{l3.1}
Let $X \ra \PP^1$ be Galois with group $G = N \rtimes B \simeq \ZZ_2^{2n} \rtimes \ZZ_3$ and action with signature $(0;\underbrace{3,\dots,3}_t)$. Then  
\begin{enumerate}
\item[(i)]  $X_N \ra \PP^1$ is totally ramified and $X \ra X_N$ is \'etale;
\item[(ii)] Over each of the $t$ branch points $p_i$ of $X \ra \PP^1$ the map $X_B \ra \PP_1$ admits $ m$ ramification points of ramification index $2$ and 
$1$ point which is \'etale over $\PP^1$.
\end{enumerate}
\end{lem}

\begin{proof}
 Let $p\in \PP^1$ be a branch point of $X \ra \PP^1$.

(i): Since $X_N \ra \PP^1$ is cyclic of degree 3, it is either \'etale over $p$
or totally ramified. Suppose it is \'etale. Since $p$ is a branch point of $X \ra \PP^1$, it is also a branch point of the Galois cover
$X \ra X_N$ of type $3$ (i.e. all points of $X$ over $p$ are of ramification degree 2 with respect to $X \ra X_N$). This gives a contradiction, since $2^{2n}$ is not divisible by 3. This completes the proof of (i).

(ii): Let $p_1,\dots, p_t$ denote the branch points of $X \ra \PP^1$ and
$b_i \in X_B$ be a point in the fibre over $p_i$. Then either  $b_i$ is of ramification degree 2 over $p_i$ or \'etale over $p_i$.
In the first case the map $X \ra X_B$ is \'etale over $b_i$ and in the second case totally ramified over $b_i$.

Let $r_i$ be the number of points in the fibre in $X_B$ over $p_i$ at which the map $X_B \ra \PP^1$ is ramified. According to the Riemann-Hurwitz \eqref{rh} formula we get
$$
g(X_B) \le \sum_{i=1}^t  r_i -(  2^{2n} -1) \le mt - (2^{2n} -1) = m(t-3)
$$ 
where the inequality follows from the inequalities $r_i \le m$. So we have
\begin{equation} \label{e3.2}
g(X_B) = m(t-3) \quad \Leftrightarrow \quad r_i = m \;\; \mbox{for all} \;\; i.
\end{equation}
On the other hand, \cite[Corollary 3.4]{r} gives a method to compute the genus of $X_B.$ Let $ \; G_j \cong {\mathbb Z}_3 \; $ be the stabilizer of the branch point $p_j \; $ with $ \; j = 1, \dots t. \; $ Then \begin{eqnarray*}
g(X_B) & = & [G : B] (0 - 1) + 1 + \frac{1}{2} \sum _{j=1}^t (\vert G : B \vert - \vert B \backslash G / G_j \vert )\\
& = & - (2^{2n} -1) +  \frac{1}{2} \sum_{i=1}^t(2^{2n} - m - 1) = m(t-3).\\
\end{eqnarray*}
Hence the equivalence \eqref{e3.2} gives the assertion.
\end{proof}

From the Riemann-Hurwitz \eqref{rh} formula we immediately get from Lemma \ref{l3.1},

\begin{lem} \label{l3.2}
Under these assumptions we have
\begin{itemize}
\item $g(X) = 2^{2n}(t -3) +1$;
\item $g(X_N) = t-2$; 
\item $g(X_{L_i}) = 2t-5$;
\item $g(X_{U_i}) = 4t-11$;
\item $g(X_{M_i}) = t-3$;
\item $g(X_B) = m(t-3)$.
\end{itemize}
\end{lem}

\begin{cor} \label{c3.3}
If $P(X_{L_i} \ra X_N)$ denotes the Prym variety of $X_{L_i} \ra X_N$, we have
$$
g(X_B) = \sum_{i=1}^m \dim(P(X_{L_i} \ra X_N) )= \sum_{i=1}^m g(X_{M_i}).
$$
\end{cor}

\begin{proof}
$\dim P(X_{L_i} \ra X_N) = g(X_{L_i}) - g(X_N)$. So Lemma \ref{l3.2} implies the assertion.
\end{proof}
Corollary \ref{c3.3}  suggests that there is a relation between these Pryms and $J(X_B)$ which we investigate in the next 
section. The second equality is explained by the following proposition.

\begin{prop} \label{p4.4}
For $i$ in $1, \ldots , m$ there is a canonical  isomorphism of principally polarized abelian varieties
$$
P(X_{L_i} \ra X_N) \stackrel{\simeq}{\lra} JX_{M_i}.
$$  
\end{prop}

\begin{proof}
Observe that for each $i$ in $1, \ldots , m$, $G/U_i \cong A_4$ acts on the corresponding curve $X_{U_i}$ (of genus $4t-11)$) and the involutions in $A_4$ act without fixed points (since $X \to X_N$ is \'etale). It then follows from the trigonal construction (see \cite{rr}) that $P(X_{L_i} \ra X_N)$ is isomorphic to $JM_i$ as principally polarized abelian varieties.  
\end{proof}

\section{The isogeny decomposition of the Jacobian of $X_B$}

Let the notations be as above.
For $i = 1, \dots,m$ denote by $\mu_i: X_B \ra X_{M_i}$ the covering given by $B \leq M_i $.  The pull back
homomorphisms 
$$
\mu_i^*: JX_{M_i} \ra JX_B
$$ 
are isogenies onto their images. Considering the composition of these isogenies with 
the addition map we get a canonical homomorphism 
$$
a: \prod_{i=1}^m JX_{M_i} \ra JX_B.
$$ 
The main result of this section is the following theorem. 

\begin{thm} \label{thm5.1}
The homomorphism $a:  \prod_{i=1}^m JX_{M_i} \ra JX_B$ is an isogeny.
\end{thm}

\begin{proof}
Recall that for any subgroup $H$ of $G$, $\rho_H$ denotes the representation of $G$ induced by the trivial representation of $H$. Moreover, $\chi_0$ denotes the trivial representation and $\theta_i,\;  i = 1, \dots,m$ the irreducible rational represetations of degree 3 of $G$. It is easy to see that
\begin{equation} \label{e5.1}
\rho_B = \chi_0 \oplus \bigoplus_{i=1}^m \theta_i. 
\end{equation}
Let $f_{B,\widetilde \chi_0}$ and $f_{B, \widetilde \theta_i}$ denote the corresponding central idempotents of the Hecke
algebra $\cH_B$. Then according to \eqref{e3.4}, the idempotent $p_B$ decomposes as 
$$
p_B = f_{B,\widetilde \chi_0} + \sum_{i=1}^m f_{B, \widetilde \theta_i}.
$$
Considering the idempotents as elements of  $\End_\QQ(JX_B)$, 
$$
JX_{B,\widetilde \chi_0} = \Ima( f_{B,\widetilde \chi_0}) = \Ima(p_G) =0,
$$ since $g(X_G) =0$ and 
$$
JX_{B,\widetilde \theta_i} = \Ima(f_{B, \widetilde \theta_i}) = \Ima(\mu_i^*),
$$
 since $\theta_i = \rho_{M_i} - \chi_0$ according to Remark \ref{r2.6}. Hence according to \eqref{e2.2},
the addition map gives an isogeny  
$$
+: \prod_{i=1}^m JX_{B,\widetilde{\theta_i}} \lra JX_B.
$$
Combining this with the isogenies $\mu_i^*$ we get the isogeny $a: \prod_{i=1}^m JX_{M_i} \ra JX_B$
as claimed.
\end{proof}

Combing several results, we can say,

\begin{cor} \label{c4.5}
For any positive integer $n \geq 2$ consider the integer $m = \frac{1}{3}(2^{2n}-1)$. Then there exist curves
of genus $m(t-3)$ for any  integer $t \ge 2n$ whose Jacobian is isogenous to the product of $m$ Jacobians.
\end{cor}

\begin{proof}
According to Lemma \ref{l4.3} there exist curves $X$ of genus $2^{2n}(t-3) +1$ with an action of the group $G = N \rtimes B$
of Section 2 and with $g(X_B) = m(t-3)$. So Theorem \ref{thm5.1} gives the assertion.
\end{proof}

\begin{cor} \label{c4.6}
Given any positive integer $N$, there exist smooth projective curves $X$ whose Jacobian is isogenous to the product of $m \ge N$ Jacobian varieties of the same dimension.
\end{cor}

\begin{proof}
Choose a positive integer $n$ such that $m = \frac{1}{3}(2^{2n} -1) \ge N$. This is equivalent to $2^{2n} > 3N$. 
According to the previous corollary there exist curves whose Jacobian is isogenous to the product of $m \geq N$ Jacobians.
 \end{proof}

\begin{rem}
The isogeny $+: \prod_{i=1}^m JX_{B,\widetilde{\theta_i}} \lra JX_B$ is the isotypical decomposition of $JX_B$ with respect to the Hecke algebra action of $\cH_B$ on $JX_B$.
\end{rem}

In order to describe the isotypical components $ JX_{B,\widetilde{\theta_k}} = \mu_k^*(JX_{M_k})$ of $JX_B$ by equations we will apply
Theorem \ref{thm3.1}.
There are $m$ conjugacy classes of involutions in $G$, and $m+1$ double cosets in $B \backslash G/ B$,
$$
G = H_0 \cup H_1 \cup \cdots \cup H_m
$$ 
with $H_0 = B$. All other $H_i$ have representatives $j_i$ in exactly one of the conjugacy classes of involutions. So 
$$
H_i = B j_i B \quad \mbox{for} \quad i = 1, \dots,m
$$
where $j_i$ is an involution of $G$.
According to Section 3, a basis of  $\mathbb{Q}[B\backslash G/ B]$ given by 
$$
q_i = \frac{1}{3} \sum_{h \in H_i} h  \quad \mbox{for} \quad  i = 0 , \dots, m.
$$ 

\begin{lem} \label{l5.3}
\begin{enumerate}
\item[(i)] $|H_i| = 9$ for $i = 1, \dots,m$;
\item[(ii)] for all $k = 1, \dots, m$ we have for all $h_i \in H_i$ and $i = 1,\dots,m$,
$$
\chi_{\theta_k}(q_0) = 1 \quad  \mbox{and} \quad  \chi_{\theta_k}(q_i) = \chi_{\theta_k}(j_i).
$$
\end{enumerate}
\end{lem}

\begin{proof} (i): Clearly $|H_i| \le 9$. Hence $3+9m = 2^{2n}3 = |G| = \sum_{i=0}^m |H_i| \le 3 + 9m$
which implies the assertion.

(ii): The first equation is obvious. For the second observe that 
$\chi_{\theta_k}(g) = 0$ for all elements $g$ of order three in $G$, and that each $H_i$ is composed by the three elements in the conjugacy class of $j_i$ and six elements of order three. This implies the assertion.
\end{proof}

The following theorem gives the description of the isotypical components of $JX_B$.

\begin{thm}
$$
\mu_k^*(JX_{M_k}) =  \{ z \in JX_B : q_i(z) = \chi_{\theta_k}(j_i)z \;\; for  \;\;  i= 1, \dots, m \}_0.
$$
\end{thm}

\begin{proof}
According to \eqref{e5.1} we have $\dim \theta_k^B = 1$. Since moreover $\theta_k$ is defined over $\QQ$, we can apply Theorem \ref{thm3.1} to give
\begin{eqnarray*}
\mu_k^*(JX_{M_k}) =  JX_{B,\widetilde{\theta_k}} &=&  \{ z \in JX :   q_i(z) = \chi_{\theta_k}(q_i)z \;\; for \;\; i= 0, \dots, m \}_0\\
&=& \{ z \in JX_B :   q_i(z) = \chi_{\theta_k}(j_i)z \;\; for \;\; i= 1,\dots,m \}_0
\end{eqnarray*}
where the last equation follows from Lemma \ref{l5.3}(ii).
\end{proof}

\section{The degree of the isogeny $a$}

In Section 5 we used the right hand side of the diagram \eqref{d3.1} to decompose the Jacobian of the curve $X_B$. By Proposition \ref{p4.4} we know that for all $i$ there is a canonical isomorphism of principally polarized abelian varieties $P(X_{L_i} \ra X_N) \ra JX_{M_i}$.  In this section we will see that one can also use the
left hand side of the diagram to decompose the Jacobian $JX_B$, here as a product of Prym varieties. In view of 
Proposition \ref{p4.4}, this is the same decomposition as above, however has the advantage that it gives in addition 
something about the degree of the isogeny.\\

Let the notation be as in the last section. For $n = 1, \dots, m$ we denote
$$
\nu_i:= X \ra X_{L_i} \qquad \mbox{and} \qquad \mu: X \ra X_B
$$
the maps of diagram \eqref{d3.1}. The maps $\nu_i^*: JX_{L_i} \ra JX$ and $\Nm \mu: JX \ra JX_B$ are the induced 
homomorphisms of the associated Jacobians. Then the addition map gives a homomorphism
$$
\alpha: \sum_{i=1}^m \Nm \mu \circ \nu_i^*: \prod_{i=1}^m P(X_{L_i} \ra X_N) \lra JX_B.
$$
According to Corollary \ref{c3.3},  $\prod_{i=1}^m P(X_{L_i} \ra X_N)$ and $JX_B$ are of the same dimension. To be more precise, we have

\begin{thm} \label{thm4.1}
$\alpha: \prod_{i=1}^m P(X_{L_i} \ra X_N) \lra JX_B$ is an isogeny with kernel contained the the $2^{2n-1}$-division points.
\end{thm}

Applying the isomorphism of Proposition \ref{p4.4}, the homomorphism  $\alpha$ coincides with the isogeny $a$ of Theorem
\ref{thm5.1}. The new result is the estimate of its degree, namely
$$
\deg(\alpha) \le 2^{2g(X_B)},
$$
which follows from the fact that the group of 2-division points of $JX_B$ is of order $2^{2g(X_B)}$.
The proof is analogous to the proof of \cite[Theorem 3.1]{clr1} (also see \cite[Proposition 3.2]{clr}).
For the convenience of the reader we give full details. We need the following proposition.

\begin{prop} \label{p4.2}
Let $f: X \ra X_N := X/N$ be a Galois cover of smooth projective curves with Galois group
$N$  and $H \subset N$ a subgroup. Denote by $\nu: X \ra Y := X/H$ and
$\varphi: Y \ra X_N$  the corresponding covers. If $\{g_1, \dots, g_r\}$ is a complete set of
representatives of $N/H$, then we have
$$
\nu^*(P(Y \ra X_N)) = \{ z \in JX^H \;|\; \sum_{i=1}^r g_i(z) = 0 \}_0.
$$
\end{prop}

\begin{proof}
This is a special case of \cite[Corollary 3.3]{clr1}
\end{proof}

Denote for $i = 1, \dots, m$,
$$
A_i := \nu_i^*( P(X_{L_i} \ra X_N)) \subset JX
$$
and let
$$
A := \sum_{i=1}^m A_i \subset JX \quad \mbox{and} \quad \widetilde A:= \mu^*(JX_B) \subset JX.
$$

Recall that $\sigma$ generates the group $B$. Then we have the following commutative diagram,
\begin{equation} \label{d2.5}
\xymatrix{
&A   \ar[rr]^{\sum_{i=0}^{2}  \sigma^i}   \ar[dr]_{\Nm \mu} && \widetilde A \ar[rr]^{\sum_{i=1}^m \sum_{h \in L_i} h}   \ar[dr]_{\beta} && A\\
\prod_{i=1}^m  P(X_{L_i} \ra X_N)   \ar[rr]_{\alpha} \ar[ur]^{\sum_{i=1}^m \nu_i^*} && JX_B \ar[rr]_{\beta \circ \mu^*} \ar[ur]_{\mu^*} && \prod_{i=1}^m P(X_{L_i} \ra X_N) \ar[ur]_{\sum_{i=1}^m \nu_i^*}
}
\end{equation}
with $\beta = (\Nm \nu_1, \Nm \nu_2, \dots, \Nm \nu_m)$.

For $i = 1, \dots,m$ consider the following subdiagram
\begin{equation} \label{d2.4}
\xymatrix{
&A_i \ar[rr]^{\sum_{i=0}^{2}  \sigma^i}  \ar[dr]_{\Nm \mu} && \widetilde A_i \ar[rr]^{\sum_{h \in L_i} h} \ar[dr]_{\Nm \nu_i} && A_i\\
 P(X_{L_i} \ra X_N)   \ar[rr]_{\alpha_i} \ar[ur]^{ \nu_i^*} && C_i \ar[rr]_{\Nm \nu_i \circ \mu^*} \ar[ur]_{\mu^*} &&  P(X_{L_i} \ra X_N) \ar[ur]_{\nu_i^*}
}
\end{equation}
with $\alpha_i:= \Nm \mu \circ \nu_i^*, \; C_i := \Nm \mu(A_i)$ and $\widetilde A_i := \mu^*(C_i)$.

\begin{prop}  \label{prop3.3}
For $i = 1, \dots, m$ the map 
$$
\Nm \nu_i \circ \mu^* \circ \alpha_i :  P(X_{L_i} \ra X_N)  \ra  P(X_{L_i} \ra X_N) 
$$
is multiplication by $2^{2n-1}$.
\end{prop}

For the proof of the proposition we need the following lemma.
\begin{lem} \label{l2.6}
For any two subgroups $H_1 \neq H_2$ of $N$  of index $2$ we have
$$|H_1 \cap H_2| = \frac{1}{2}|H_1| = 2^{2n-2}.$$
\end{lem}

\begin{proof}
Since $N = H_1H_2$, we have
$N/H_1 \simeq H_2/(H_1 \cap H_2),$
which implies the assertion.
\end{proof}

\begin{proof}[Proof of Proposition \ref{prop3.3}]
Since $\nu_i^*:  P(X_{L_i} \ra X_N)  \ra A_i$ is an isogeny, it suffices to show that the composition
$$
\Phi_i: = \sum_{h \in L_i} h \circ \sum_{i = 0} ^{2} \sigma^i: A_i \ra A_i
$$
is multiplication by $2^{2n-1}$.

Now from Proposition \ref{p4.2} we deduce
\begin{equation} \label{e2.7}
A_i = \{ z \in JX \;|\; hz = z \; \mbox{for all} \; h \in L_i \; \mbox{and} \; nz = -z \; \mbox{for all} \; n \in N \setminus L_i  \}_0
\end{equation}
For the proof of \eqref{e2.7} we apply Proposition \ref{p4.2} with $H = L_i$. Let $g_1$ be any element of $L_i$
and $g_2$ be any element of $N \setminus L_i$; then $\{ g_1 , g_2\}$ is a complete set of representatives of $N/L_i$. So Proposition \ref{p4.2} gives
$$
A_i = \{ z \in JX^{L_i} \; |\; g_1(z) + g_2(z) = 0 \}_0.
$$ 
But $z \in JX^{L_i}$ implies that $g_1(z) = z$. This gives
$$
A_i = \{ z \in JX^{L_i} \; |\; g_2(z) = -z \}_0.
$$
But $X_{L_i}/X_N$ is a double cover, which implies $g_2(z) = -z$. Since $g_1$ and $g_2$ are arbitrary elements of $L_i$ and $N \setminus L_i$ respectively, this gives the assertion.

Now for any $z \in A_i$,
$$
\Phi_i(z) = \sum_{h \in L_i} h(z) + \sum_{h \in L_i} h \sum_{k=1}^{2} \sigma^k(z).
$$
By equation \eqref{e2.7} and Lemma \ref{l2.6} we have
$$
\sum_{h \in L_i} h (z) = |L_i| z = 2^{2n-1}z
$$
and for $k = 1$ and 2,
$$
\sum_{h \in L_i} h \sigma^k(z) = \sigma^k \sum_{h \in \sigma^{-k} L_i \sigma^k} h(z) = 0,
$$
since $L_i \neq \sigma^{-k} L_i \sigma^k$ and by Lemma \ref{l2.6}  half of the elements of the
subgroup $\sigma^{-k} L_i \sigma^k$ belong to $L_i$, hence fix $z$, and the other half belongs
to $N \setminus L_i$ and hence sends $z$ to $-z$. Together they complete the proof of the proposition.
\end{proof}

\begin{proof}[Proof of Theorem \ref{thm4.1}]
Since
$$
\beta \circ \mu^* \circ \alpha = \prod_{i=1}^m (\Nm_{\nu_i} \circ \mu^* \circ \alpha_i),
$$
Proposition \ref{prop3.3} implies that $\beta \circ \mu^* \circ \alpha$ is multiplication by $2^{2n-1}$.
In particular $\alpha$ has finite kernel. But according to Corollary \ref{c3.3},
$\prod_{i=1}^m  P(X_{L_i} \ra X_N)$ and $JX_B$ have the same dimension. So $\alpha$ is an isogeny.
\end{proof}

\begin{rem}
As we noticed already above, the isogenies $a$ and $\alpha$ are compatible with respect to the isomorphism of
Proposition \ref{p4.4}. Hence Theorem \ref{thm4.1} is also the isotypical decomposition of $JX_B$.
\end{rem}

\section{The group algebra decomposition of $JX$}

Let $G = N \rtimes B$ as act on the curve $X$ as above. 
For the trivial subgroup $H = \{1\}$ the Hecke algebra $\cH_{\{1\}}$ coincides with the whole group algebra $\QQ[G]$
and we have $JX_{\{1\}} = JX$ and $\widetilde W_i = W_i$ for all $i$.  According to \eqref{e2.2} and \cite[equation (2.3)]{clr1}
the isotypical decomposition of $JX$ is
$$
JX \sim JX_\psi \times \prod_{i=1}^m JX_{\theta_i},
$$
where we write for simplicity $JX_{\{1\},W} = JX_{W}$ for each irreducible rational representation $W$ of $G$, for which we refer to Corollary \ref{c2.4}.  Note that $JX_{\chi_0} = 0$, since $X_G = \PP^1$.

This decomposition can be decomposed further. In fact, as outlined in \cite[Section 2]{clr1}, for each $\theta_i$ there is a 
not uniquely determined abelian subvariety $B_i$ such that $JX_{\theta_i} \sim B_i^3$. On the other hand, at least by these 
methods $JX_\psi$ cannot be decomposed further. Hence we get a decomposition
$$
JX \sim JX_\psi \times \prod_{i=1}^m B_i^3,
$$
which is called the {\it group algebra decomposition} of $JX$ (see \cite[Section 13.6]{bl}). 

\begin{lem} \label{l7.1}
\begin{enumerate}
\item[(i)] $JX_\psi = JX_N$;
\item[(ii)] $B_i \sim P(X_{L_i} \ra X_N)$ for $i = 1, \dots, m$.
\end{enumerate}
\end{lem}

\begin{cor}
A group algebra decomposition of $JX$ is 
\begin{eqnarray*}
JX &\sim & JX_N \times \prod_{i=1}^m (P(X_{L_i} \ra X_N))^3.\\
&\sim & JX_N \times \prod_{i=1}^m (JX_{M_i})^3
\end{eqnarray*}
\end{cor}

\begin{proof}[Proof of the Lemma \ref{l7.1}] (i): The covering $X_N \ra \PP^1$ is Galois with group of order 3. There are 2
irreducible rational representation of this group, namely $\psi$ and the trivial representation $\chi_0$. To $\psi$ corresponds the
Prym variety $P(X_N \ra \PP^1)$ and to the trivial representation $J\PP^1 =0$. This implies $JX_\psi \sim JX_N$ and thus  
$JX_\psi = JX_N$, since the isotypical decomposition is uniquely determined. 

(ii): According to Theorem \ref{thm4.1}, $\alpha(P(X_{L_i} \ra X_N))$ contained in $JX_B$ and thus in $JX$, since $X \ra X_B$ 
is ramified. So clearly it is contained in $JX_{\theta_i}$. Since one can conclude from Lemma \ref{l3.2} that 
$\dim JX_{\theta_i} = \dim B_i$, this implies that there is an isogeny $JX_{\theta_i} \sim B_i$.

The last isogeny is a consequence of  Proposition \ref{p4.4}.
\end{proof}

\end{document}